\documentclass{article}

\def\P{{\mathcal P}}

\def\S{{\mathcal S}}
\def\V{{\mathcal V}}

\newfont{\bb}{msbm10 at 10pt}
\def\r{\hbox{\bb R}}

\def\z{\hbox{\bb Z}}
\def\x{\hbox{\bf x}}

\usepackage[latin1]{inputenc}

\newenvironment{proof}{\trivlist
\item[\hskip\labelsep{\it Proof}\,:]}{\hfill{$q.e.d.$}\endtrivlist}
\newtheorem{theorem}{Theorem}[section]
\newtheorem{lemma}[theorem]{Lemma}

\newtheorem{corollary}[theorem]{Corollary}
\newtheorem{definition}[theorem]{Definition}
\newtheorem{remark}[theorem]{Remark}
\begin{document}

\title{Capillary channels in a gravitational field\footnote{Partially
supported by MEC-FEDER 
 grant no. MTM2004-00109.}}
\author{Rafael L\'opez\\
Departmento de Geometría y Topología\\
Universidad de Granada\\
18071 Granada (Spain)\\
e-mail:{\tt rcamino@ugr.es}}
\date{}

\maketitle

\begin{abstract}
The liquid shape between two vertical parallel plates in a gravity field due to 
capillary forces is studied. When the physical system  achieves its
mechanical equilibrium, the capillary surface has mean curvature proportional to its height
above a horizontal reference plane and it meets the vertical walls in a
prescribed angle.  We examine the shapes of these interfaces and their qualitative properties depending on the sign of the capillary constant. We focus to obtain estimates of the size 
of the meniscus, as for example, its height and volume.
\end{abstract}

%%%%%%%%%%%%%%%%%%%%%%%%%%%%%%%%
\section{Introduction and formulation of the problem}
%%%%%%%%%%%%%%%%%%%%%%

 Consider an infinite horizontal reservoir of fluid and 
let us introduce two  vertical parallel plates. 
The action of capillarity causes that the liquid rises between both plates
until a state of mechanical equilibrium. Denote $\S$ the interface liquid-air formed 
by the fluid between the two plates and whose shape we would like to determined.  The fluid surface level at large distance from the plates provides a reference level  
$\Pi$ for atmospheric pressure that does not change with perturbations of the fluid surface between the plates. 
According to the principle of virtual work,  the  configurations that adopts the liquid between 
the two plates are characterized by two facts \cite{fi1}: 
\begin{enumerate}  
\item The mean curvature of  $\S$ is proportional to 
the height of $\S$  with respect to $\Pi$ (Laplace  equation).
\item The  angles $\gamma_i$ with which $\S$ intersects 
the plates  are constant (Young condition). These constants depend only on the materials of the liquid and the plates.
\end{enumerate}

We also can consider the wetting phenomenon when one spreads out a sufficient amount of liquid on a stripped domain in such way the  liquid tends to wet the domain.

Consider $(x,y,z)$ the usual coordinates in Euclidean three-space
$\r^3$,  $P_i=\{x=a_i\}$ the two vertical planes with $a=a_1=-a_2>0$ and    
 $\Pi=\{z=0\}$ the horizontal plane. Set  $L_i=\Pi\cap P_i$. Denote $\Omega=\{(x,y)\in\r^2;|x|<a\}$ the horizontal strip 
in $\Pi$ determined by the two  planes, identifying $\r^2$ with $\Pi$ as usually.
Let the height of this capillary free surface $\S$ with respect to $\Pi$, 
assumed nonparametric over  $\Omega$, be given by the
scalar function $u=u(x,y)$, $(x,y)\in\Omega$. When the capillary and gravity forces are in 
equilibrium, 
$u$ satisfies the partial differential equation
\begin{equation}\label{capi1}
\mbox{div }Tu=\kappa u
\end{equation}
in $\Omega$ where
$$Tu=\frac{Du}{\sqrt{1+|Du|^2}}.$$
See \cite{fi1}. Here  $\kappa=\rho g/\sigma$ is the capillarity constant with $\sigma$, the surface tension, $\rho$, the difference of densities across the interface $\S$ and $g$, the gravitational acceleration, with positive and negative sign in the 
sessile and pendent case respectively. Equation (\ref{capi1}) can be interpreted as that   the mean curvature $H$ of the surface $z=u(x,y)$ is 
$\kappa u/2$. The  Young condition writes as
\begin{equation}\label{capi2}
\nu_i\cdot Tu=\cos\gamma_i\hspace*{1cm}\mbox{along $L_i$}
\end{equation}
where $\nu_i$ is the  unit exterior normal on $L_i$. Here $\gamma_i$ are the  contact angles
  with which $\S$
meets  $P_i$, $i=1,2$. The orientation on $\S$  points in the $z$-positive direction. 
If the two plates are made with the same materials, $\gamma=\gamma_1=\gamma_2$.
We may normalize so  $0\leq\gamma\leq \pi$. The range $0\leq\gamma\leq\pi/2$  indicates capillary rise; $\pi/2<\gamma\leq\pi$ yields
capillary fall. The angles $\gamma_i$ are determined by the  volume per unit of length
enclosed by the surface. If $\Omega_b=(-a,a)\times (-\frac{b}{2},
\frac{b}{2})$ is a rectangular piece of $\Omega$ of length 
$b$, an integration of (\ref{capi1}) gives
$$\kappa V_b=(\cos\gamma_1+\cos\gamma_2)b+4a,$$
where $V_b$ is the enclosed volume by $u$ over the domain $\Omega_b$. Then 
$$\cos\gamma_1+\cos\gamma_2=\lim_{b\rightarrow\infty}\kappa\frac{V_b}{b}.$$
This identity is similar  when $\Omega$ is a bounded domain, namely, 
$$\cos\gamma=\kappa\frac{V_{\Omega}}{|\partial\Omega|}.$$
When the effect of gravity is ignored, the liquid-air interface is characterized by 
a constant mean curvature 
surface. A first example of graph with constant mean curvature on a band is
 any section of an infinite round cylinder positioned with its axis parallel to $P_i$. 
Exactly, this example motivates us to consider that the shape of the surface $\S$ is translationally invariant with respect to 
the $y$-coordinate. Our surface will be invariant by the reflection with  respect to the plane $\{y=t\}$ and $\S$ is determined by its intersection with any plane $\{y=t\}$. 
Then $\S$ is a cylindrical ruled surface.

Classically,  it has been  studied the   capillary problem when
 the liquid rises in a tube with circular section  in such way that our setting reduces then to consider one of the curvature radius is infinite.  So, in the literature, the capillary problem studied here it has been considered 
  in the study of the shape of a meniscus facing a vertical plate.   
As we will see, in the one-dimensional problem,
a first integration of (\ref{capi1}) is 
obtained in such way that the solutions can be expressed in terms of elliptic 
integrals and some estimates of the height of the meniscus have been  obtained from these integrals or as limit case of
 the two-dimensional case \cite{ba,bak,gbq,he,ne,si,vo}. Our interest in this work is twofold. First, we analyze the 
symmetries of the surface and the shapes adopted  depending on the sign of $\kappa$. On the 
other hand, we shall obtain a detailed study of the height of the meniscus, as well as, other estimates on the volume. In this sense, we will follow the same spirit as in \cite{cf1,fi0,fi2,fi3}. See also \cite{fi1}. 

This paper is organized as follows. 
In Section \ref{planes} we describe all cylindrical ruled surfaces whose mean curvature is proportional to its height 
with respect to $\Pi$, making a study of their symmetries. In Section \ref{estimates1}, we  obtain estimates of the height of the meniscus in the capillary problem. In 
Section \ref{estimates2} we consider  sessile  liquid channels 
with results on existence with respect to the volume enclosed by the  channel. Finally in Section  \ref{estimates3} 
we study pendent liquid channels with the main   conclusion 
 that the morphologies that will appear differ  completely  than in the two-dimensional problem.

%%%%%%%%%%%%%%%%%%%%%%%%%%%%%%%%%%%%%%%%%%
\section{Capillary immersed bands}\label{planes}
%%%%%%%%%%%%%%%%%%%%%%%%%%%%%%%%%%%%%%%%%%%

With  more of generality, let $S$ be a cylindrical ruled surface in the space $\r^3$, that is, an immersed surface
parametrized as 
$$\x(s,t)=\alpha(s)+t \vec{w},\hspace*{1cm}s\in I, t\in\r,$$
where $\alpha$ is a regular planar  curve of $\r^3$ defined in some interval $I$, called the directrix of $S$ and  $\vec{w}\in\r^3$, $|\vec{w}|=1$. 
In this section we are interested by those cylindrical surfaces that satisfy the capillary equation (\ref{capi1}) for 
some constant $\kappa\not=0$.
Without loss of generality, we assume that $\alpha$ is  parametrized by arc length,   $\langle\alpha',\vec{w}\rangle=0$ and  the binormal of $\alpha$ in the Frenet trihedron is $-\vec{w}$. 
Then  $H=C_{\alpha}/2$, where $C_{\alpha}$ is the curvature of $\alpha$. 
Equation (\ref{capi1}) implies then 
$$C_{\alpha}(s)=\kappa z(\alpha(s))+t z(\vec{w})$$ 
for all $s$ and $t$. We infer  that $z(\vec{w})=0$, that is, $\vec{w}$ is a horizontal vector and the rulings of the 
surface are horizontal straight-lines. 

\begin{definition} Let $\kappa\not=0$. A $\kappa$-cylindrical surface is a cylindrical ruled surface that 
locally satisfies the capillary equation (\ref{capi1}).
\end{definition}

 In particular, 
each vertical plane orthogonal to the rulings is a plane of symmetry of $S$. In addition,   the angle that makes a such surface with a pair of vertical parallel planes, or  
a horizontal plane, is constant. For the study of  the existence of $\kappa$-cylindrical surfaces, 
we parametrize the surface $S$ as  $\x(s,t)=(x(s),t,z(s))$,  $s\in I, t\in \r$,
where $\alpha(s)=(x(s),z(s))$. We know that
\begin{equation}\label{unit}
x'(s)^2+z'(s)^2=1,\hspace*{1cm}s\in I.
\end{equation}
Let $\theta(s)$ be the angle between the vectors $\partial/\partial x$ and 
 $\alpha'(s)$. By  (\ref{unit}), the equation (\ref{capi1}) converts into 
 the O.D.E. system $\P$:
\begin{eqnarray}
x'(s)&=&\cos\theta(s)\label{eq1}\\
z'(s)&=&\sin\theta(s)\label{eq2}\\
\theta'(s)&=&\kappa z(s)\label{eq3}
\end{eqnarray}

\begin{theorem}\label{exis} The system of ordinary differential equations $\P$ has a unique solution for each initial condition. Moreover the maximal interval of 
the solution is $\r$.
\end{theorem}

\begin{proof}
Classical theory yields existence of solutions for each initial data $x(0)=x_0, z(0)=z_0,
\theta(0)=\theta_0$. Denote $\P(x_0,z_0,\theta_0)$ the 
initial value problem for the initial conditions $(x_0,z_0,\theta_0)$. If $(x,z,\theta)$ is a solution for 
$\P(x_0,z_0,\theta_0)$, then $(x+a,z,\theta+b)$ is the solution of $\P(x_0+a,z_0,\theta_0+b)$. Thus, we can assume
that $x_0=0$ and $\theta_0=0$. Let us denote $\P=\P(z_0)$. In this article, we assume  these initial conditions.

For $z(0)=z_0$, we obtain a solution $(x,z,\theta)$ defined around $s=0$. 
It is immediate that $(s,0,0)$ is the solution for $z_0=0$. Assume now $z_0\not=0$. From 
(\ref{eq1})-(\ref{eq3}), 
$$x''=-\theta' z'=-\kappa z z'=-\frac{\kappa}{2}(z^2)'.$$
Then there exists a constant $m\in\r$ such that 
$$x'=-\frac{\kappa}{2}z^2+m=\cos\theta.$$
At $s=0$, we have $m=1+\kappa z_0^2/2$. Thus
$$\frac{\kappa}{2}z^2=1-\cos\theta+\frac{\kappa}{2}z_0^2$$
or
\begin{equation}\label{cotaz}
z(s)^2=z_0^2+\frac{2}{\kappa}\left(1-\cos\theta(s)\right).
\end{equation}
Therefore $z$ is a bounded function. As a consequence of (\ref{cotaz}), together
with  (\ref{eq1})-(\ref{eq3}), 
 the first derivatives of $x$, $z$ and $\theta$ are bounded 
functions and
the theory of O.D.E. yields that the solutions can be continued indefinitely. 
This proves the result. 
\end{proof}

We prove that our  $\kappa$-cylindrical surfaces 
have a rich symmetry. 

\begin{theorem}[Symmetry I] \label{t-i}
Let  $S\subset\r^3$ be a $\kappa$-cylindrical surface. Then $S$ is symmetric with respect to 
 any vertical plane parallel to the rulings and 
that acrosses an 
extremum of the height of the function $z$, where $\alpha=(x,z)$.
\end{theorem}
 
\begin{proof} 
Consider $\alpha$ the directrix curve of $S$, and we assume
the initial data given in Theorem \ref{exis}. It suffices to  prove that the trace of $\alpha$ is symmetric with 
respect to the line $x=x(s_0)$, where $s_0$ is any value with $\cos\theta(s_0)=\pm 1$. Let $m\in\z$ be such 
that $\theta(s_0)=m\pi$. The theorem is proved if  for $s\in\r$,
\begin{eqnarray*}
x(s+s_0)-x(s_0)&=&x(s_0)-x(s_0-s)\\
 z(s+s_0)&=&z(s_0-s)\\
 \theta(s+s_0)&=&2m\pi-\theta(s_0-s).
\end{eqnarray*}
However, these two sets of functions satisfy the 
same O.D.E. system $\P$ and  initial conditions. The uniqueness of solutions of an O.D.E.  concludes  the proof.

\end{proof}

In a similar way, we have

\begin{theorem}[Symmetry II] \label{t-ii}
Let  $S\subset\r^3$ be a $\kappa$-cylindrical surface and 
$\alpha=(x,z)$ its directrix. Assume that $\alpha$
accrosses the $x$-axis at $s=s_0$. Then $\alpha$ is symmetric with respect to 
the point $\alpha(s_0)$. 
\end{theorem}

To end this section, we study the symmetries of the surfaces by 
horizontal translations orthogonal 
to the rulings. We distinguish two cases depending on the sign of $\kappa$.

\begin{theorem}[Sessile case]\label{t-sessile} Let $S$ be a $\kappa$-cylindrical surface with $\kappa>0$.
Then there exists a horizontal vector $\vec{v}$ orthogonal to the rulings such that $S$
is invariant by the group of translations 
$$G=\{t_n;n\in\z\},\hspace*{1cm}t_n(p)=p+n\vec{v}.$$
Moreover, if $\alpha=(x,z)$ is the directrix of $S$, the function $z$  is a periodic.
\end{theorem}

\begin{proof}  Since that $(x,-z,-\theta)$ is a solution 
of $\P(-z_0)$ provided $(x,z,\theta)$ is the one of the system $\P(z_0)$, we assume $z_0>0$.
From (\ref{cotaz}), $z\geq z_0$. Equation (\ref{eq3}) implies
that $\theta$ is  strictly increasing and its limit is $\infty$. Set $T>0$ the first  number such that 
$\theta(T)=2\pi$. Again, the uniqueness of solutions in a O.D.E. gives 
$\alpha(s+T)=\alpha(s)+(x(T),0)$. 
 This  means that  the  surface is invariant by 
the group of translations $G$, with  $\vec{v}=
(x(T),0,0)$.

\end{proof}

\begin{remark}\label{rem1}
As consequence of Theorem \ref{t-sessile}, and because 
 $\theta$ is increasing function to infinity, the 
velocity vector rotates infinite times around the origin. 
\end{remark}

From (\ref{cotaz}) and because $\cos\theta$ takes all the values into
the interval $[-1,1]$, we bound the height function $z$ in terms of the 
 lowest height $z_0$.

\begin{corollary} Let $S$  be a $\kappa$-cylindrical surface with $\kappa>0$ and denote $z$ the height
with respect to the plane $\Pi$. Then $z$ satisfies
$$z_0\leq z(p)\leq \sqrt{\frac{4}{\kappa}+z_0^2},\hspace*{1cm}p\in S,$$
and both  bounds are achieved.
\end{corollary}

This means that, fixed $\kappa$, the difference between the two extremum of $z$ 
is bounded by a  constant, namely $\sqrt{4/\kappa}$, independent on $z_0$.

\begin{theorem}[Pendent case] \label{p-sessile} Let $S$ be a $\kappa$-cylindrical surface
with $\kappa<0$. Denote by  $\alpha=(x,z)$ its  directrix obtained as solution of $\P(z_0)$. 
Without loss of generality, suppose $z_0<0$. 
\begin{enumerate}
\item If $z_0<-2/\sqrt{-\kappa}$, then there exists horizontal 
vector $\vec{v}$ orthogonal to the rulings such that $S$
is invariant by the group of translations 
$G=\{t_n;n\in\z\}$, with $t_n(p)=p+n\vec{v}$. 
Moreover $z$ is a periodic function and $z<0$. 
\item If $z_0=-2/\sqrt{-\kappa}$, then $z_0\leq z<0$, $z$ is strictly increasing
and $\lim_{s\rightarrow\infty}z(s)=0$.
\item If $-2/\sqrt{-\kappa}<z_0<0$, then $\alpha$ is a periodic function.
Moreover, $z$ vanishes in a discrete set of 
point, $z_0\leq z(s)\leq -z_0$  where both extremum are achieved and $\alpha$ is symmetric with respect to any zero of $z$.
\end{enumerate}
\end{theorem}

\begin{proof} 
\begin{enumerate}
\item  Identity (\ref{cotaz}) implies
 that $z$ does not vanish and 
\begin{equation}\label{zeta1}
z_0\leq z(s)\leq -\sqrt{z_0^2+\frac{4}{\kappa}}.
\end{equation}
 From (\ref{eq3}), $\theta$ is strictly increasing with 
$$\theta'\geq-\kappa \sqrt{z_0^2+\frac{4}{\kappa}}.$$
This means that $\theta$ increases until infinity. Again, let $T>0$ be the first number where $\theta(T)=2\pi$. The 
same reasoning as in Theorem \ref{t-sessile}, proves the statement 1. In particular, $\cos\theta(s)$ takes all values in 
$[-1,1]$ and the bounds in (\ref{zeta1})  are achieved.

\item From (\ref{cotaz}),  the  only zeroes of $z$  occur when $\cos\theta(s)=-1$.
If $z$ vanishes at some point, the uniqueness of solutions would imply that $z=0$, which is a contradiction. Thus 
$z<0$ and $\cos\theta>-1$. Near $s=0$, $\theta$ is increasing and the same occurs for 
the function $z$. In addition, $0\leq\theta(s)<\pi$. 
Moreover, $z(s)<0$,  $z'(s)>0$ for $s\in \r$ and 
$$\lim_{s\rightarrow \infty}z(s)=z_1\hspace*{.5cm}\lim_{s\rightarrow \infty} z'(s)=0,$$ 
for some number $z_1\leq 0$. If $z_1<0$, by (\ref{eq3}) $\theta'>k>0$, for some constant $k$ and 
$\theta$ would attain the value $\pi$. This contradiction 
yields $z_1=0$. 

\item First, we prove that $z$ vanishes. Because $z_0<0$, the functions 
$z$ and $\theta$ are  increasing near $s=0$. If $z\leq 0$,  then 
$$\lim_{s\rightarrow\infty}z(s)=\delta\hspace*{1cm}
\lim_{s\rightarrow\infty}z'(s)=0,$$
for some  number $\delta\leq 0$.
If $\delta<0$, (\ref{eq3}) implies that $\theta$ increases until $\infty$ 
and $\cos\theta$ takes all possible value. Equation  (\ref{cotaz}) 
together with $-2/\sqrt{-\kappa}<z_0$ imply that $z=0$ at some point. If
$\delta=0$, (\ref{eq3}) gives again that either $\theta\rightarrow\infty$, which is 
a contradiction or $\theta\rightarrow\theta_0$, for some $\theta_0<\infty$.
Letting $s\rightarrow \infty$ in (\ref{eq2}), we conclude that $\theta=\pi$, 
in contradiction with (\ref{cotaz}) and $-2/\sqrt{-\kappa}<z_0$.

Therefore, $z$ vanishes at some point. We use the Theorem \ref{t-ii} to conclude 
that $z$ is symmetric with respect to any zero of $z$. Furthermore, 
$z$ has a minimum at $s=0$, since 
$z''(0)=\kappa z_0>0$. Then $z$ is a bounded function with
 $z_0\leq z\leq -z_0$. Moreover, the same Theorem \ref{t-ii} yields $z(2s_0)=-z_0$. 
Then (\ref{cotaz}) implies that at $s=2s_0$ (resp. $s=0$), $z$ attains a maximum (resp. minimum). The proof finishes using the symmetries of $(x(s),z(s))$ given
in Theorem \ref{t-i}.  Exactly, it follows that
\begin{eqnarray}\label{perio}
x(s+4s_0)&=&x(s)+x(4s_0)\\
z(s+4s_0)&=&z(s)\\
\theta(s+4s_0)&=&\theta(s)
\end{eqnarray}
 \end{enumerate}

\end{proof}

%%%%%%%%%%%%%%%%%%%%%%%%%%%%
\section{Estimates of capillary strips: case $\kappa>0$}\label{estimates1}
%%%%%%%%%%%%%%%%%%%%%%%%%%%%%%%

In this section, we consider $\kappa$-cylindrical surfaces $S$ that are graphs over the strip $\Omega$ of a  function $u$, that is, $S$ is the surface 
$z=u(x,y)$ that projects simply onto $\Omega$. We will derive estimates for the  capillary rise, as for example, the center height $u_0$ and the outer height $u(a)$. 
For the two dimensional problem, we refer \cite{fi2,fi3,fi1}.

 Setting $r=x$ and 
$u(r,y)=u(r)$, Equation (\ref{capi1}) becomes 
\begin{equation}\label{capi3}
\frac{u''(r)}{(1+u'(r)^2)^{3/2}}=\frac{d}{dr}\left(\frac{u'(r)}{\sqrt{1+u'(r)^2}}\right)=\kappa u.
\end{equation}
Together  (\ref{capi3}), we consider  the initial conditions 
\begin{equation}\label{capi33}
u(0)=u_0>0,\hspace*{1cm}u'(0)=0.
\end{equation}
Denote $u=u(r;u_0)$ the dependence of $u$ with respect to the initial condition  $u(0)=u_0$.
It is immediate then that 
\begin{enumerate}
\item  $u(r;0)=0$ and $u(-r;u_0)=u(r;u_0)$.
\item $u(r;u_0)=-u(r;-u_0)$: up a symmetry with respect to the $r$-axis, 
the sign of the initial condition $u_0$ can be prescribed.
\end{enumerate}

According to these properties, we will assume that $u_0\not=0$ and
that the signs of $u_0$ and $\kappa$ agree. Although much of our results are valid
with independence on the sign of $\kappa$, we restrict
to the case  that $\kappa$ is a positive number.

The boundary condition (\ref{capi2}) writes now
$u'(a)=\cot\gamma$. We know by standard theorems of O.D.E. that there exists such function $u$ defined in some interval 
around $r=0$. Put
$$\sin\psi=\frac{u'}{\sqrt{1+u'^2}},\hspace*{1cm}\cos\psi=\frac{1}{\sqrt{1+u'^2}},$$
where $\psi$ is the  angle that makes $u(r)$ with the horizontal line at each point  $(r,u(r))$.
 Then (\ref{capi3}) takes the 
form 
\begin{equation}\label{B1}
(\sin\psi)'=\kappa u.
\end{equation}
For $r>0$ and close to $0$, 
$$\sin\psi=\kappa \int_0^r u(t)\ dt.$$
As $u_0>0$, the integrand is positive near to $r=0$. Then $\sin\psi>0$, and so, $u'(r)>0$.
This means that $u$ is increasing provided that $u$ is defined
in the maximal interval $(0,R)$, $0<R\leq\infty$. 
Multiplying by $u'$ in (\ref{capi3}), we have a first integration
$$\frac{1}{\sqrt{1+u'^2}}=-\frac{\kappa}{2}u^2+c,$$
for some constant $c$. At $r=0$, 
$$c=1+\frac{\kappa}{2}u_0^2.$$

Therefore
\begin{equation}\label{first}
u'=\sqrt{\frac{4}{\left(2+\kappa(u_0^2-u^2)\right)^2}-1}
\end{equation}
and 
$$u(r)=u_0+\int_0^r\sqrt{\frac{4}{(2+\kappa(u_0^2-u(t)^2)^2}-1}\ dt.$$
From (\ref{capi3}) and (\ref{capi33}), we have $u''\geq \kappa u\geq \kappa u_0>0$. This 
implies that $u'$ in increasing on $r$ and $u'(R)=\infty$. 
Equation (\ref{first}) (or (\ref{cotaz})) gives
$$u(R)=\sqrt{\frac{2}{\kappa}+u_0^2}.$$
This means that $R<\infty$ and that 
the maximal distance between the center and outer height
 of a $\kappa$- cylindrical surface is
\begin{equation}\label{uR}
u(R)-u_0=\frac{\frac{2}{\kappa}}{u_0+\sqrt{\frac{2}{\kappa}+u_0^2}}.
\end{equation}
This was to be expected according to  the Remark \ref{rem1} 
and (\ref{cotaz}). As a consequence, 
fixed a capillarity constant $\kappa>0$ and $u_0>0$, the 
angle of contact $\gamma$ takes all the values in the range $0\leq\gamma\leq\pi/2$.
More generally, we have from (\ref{cotaz})

\begin{corollary}
Let $S$ be a $\kappa$-cylindrical surface  given by the profile $u=u(r;u_0)$. 
If $\gamma$ is the contact angle with the vertical walls at $r=a$, then 
\begin{equation}\label{q1}
q:=u(a)-u(0)=\frac{\frac{2}{\kappa}(1-\sin\gamma)}{u_0+\sqrt{u_0^2+\frac{2}{\kappa}(1-\sin\gamma)}}<
\sqrt{\frac{2}{\kappa}(1-\sin\gamma)}.
\end{equation}
\end{corollary}

Fixing $\gamma$, the function $q=q(u_0)$ depending on the initial condition $u_0$ goes from $0$ to $\sqrt{2(1-\sin\gamma)/\kappa}$,
with
$$\lim_{u_0\rightarrow 0} q=\sqrt{\frac{2}{\kappa}(1-\sin\gamma)},
\hspace*{1cm} \lim_{u_0\rightarrow \infty} q=0.$$

As $u$ is increasing on $r$, we bound the integrand in (\ref{B1}) by $u_0<u(t)<u(r)$ obtaining 
\begin{equation}\label{B2}
 \kappa   u_0<\frac{\sin\psi(r)}{r}<\kappa  u(r).
\end{equation}
Moreover, 
$$\lim_{r\rightarrow 0}\frac{\sin\psi(r)}{r}=\kappa u_0.$$
This allows to give the following results on existence

\begin{theorem}\label{sessile-ex}  Let $\kappa>0$ be a constant of capillarity. Given $2a>0$, the width of the strip $\Omega$, and 
$0\leq\gamma<\pi/2$, a contact angle, 
there exists  a $\kappa$-cylindrical surface 
on $\Omega$ that makes a contact angle $\gamma$ with the plates $P_1\cup P_2$.
\end{theorem}

\begin{proof} 

The problem reduces to find $u_0>0$ such that $u'(a;u_0)=\cot\gamma$,
or in terms of the function $\psi$, that 
$\sin\psi(a)=\cos\gamma$, where $0<\cos\gamma\leq 1$.
If $u_0=0$, we know that $u(r;0)=0$. By the continuity on the parameter
$u_0$ for the solutions
of (\ref{capi3})-(\ref{capi33}),
$$\lim_{u_0\rightarrow 0}\sin\psi(a;u_0)=\sin\psi(a;0)=0.$$
If we denote by $R=R(u_0)$ the maximal interval of 
$u(r;u_0)$, and since $R(0)=\infty$, there exists $u_0$ close to  $0$ 
such that the following holds:
$$R(u_0)>a\hspace*{1cm}\sin\psi(a;u_0)<\cos\gamma.$$
From (\ref{B2}), $u_0$ cannot take any value, but its supremum is
$1/(\kappa a)$. Again, the left inequality in (\ref{B2}) leads to
$$\lim_{u_0\rightarrow 1/(\kappa a)}\sin\psi(a;u_0)=1.$$
By continuity, there exists $u_0\in(0,1/(\kappa a))$, such 
that the solution $u(r;u_0)$ satisfies $\sin\psi(a)=\cos\gamma$.

\end{proof}

Now, we bound the center height $u_0$ and the outher height $u(a)$. 
Consider a lower circular arc $\Sigma^{(1)}$: $u^{(1)}(r)$, centered on the $u$-axis, with 
$u^{(1)}(0)=u_0$ and of radius $R_1=1/(\kappa u_0)$.  Let also $\Sigma^{(2)}$: $u^{(2)}(r)$ 
be a lower circular arc, centered on the $u$-axis, with $u^{(2)}(0)=u_0$, and such that 
$$\frac{d}{dr}u^{(2)}(a)=\frac{d}{dr}u(a)$$
so that $\Sigma^{(2)}$ meets the vertical plates in the same angle as does the solution surface. 

The circular arcs $\Sigma^{(1)}$ and  $\Sigma^{(2)}$ can parametrized as
\begin{equation}\label{u1}
u^{(1)}(r)=u_0+R_1-\sqrt{R_1^2-r^2},\hspace{1cm}R_1=\frac{1}{\kappa u_0}.
\end{equation}
\begin{equation}\label{u2}
u^{(2)}(r)=u_0+R_2-\sqrt{R_2^2-r^2},\hspace*{1cm}R_2=\frac{a}{\cos\gamma}.
\end{equation}

{\it Claim.} The three functions satisfy $u^{(1)}(r)<u(r)<u^{(2)}(r)$ in the interval $(0,a]$.

\begin{proof}[of the Claim]
By (\ref{capi3}), the  curvature of  $u(r)$  is 
$$C_u(r)=\frac{u''(r)}{(1+u'(r)^2)^{3/2}}=\kappa u(r).$$
Moreover $C_u$ is increasing on $r$ since $\kappa$ and $u'$ are positive. At $r=0$, 
$C_u(0)=\kappa u_0=C_{u^{(1)}}(0)$ and $\Sigma^{(1)}$ has constant curvature. Because
$u(0)=u^{(1)}(0)$, we conclude then
$$\frac{d}{dr}u^{(1)}(r)<\frac{d}{dr}u(r);\hspace*{.5cm}u^{(1)}(r)<u(r);\hspace*{.5cm}0<r<a.$$
For $u^{(2)}$, as $u^{(1)}$ and $u^{(2)}$ are circles and  $C_{u^{(2)}}(a)>C_{u^{(1)}}(a)$, then 
$C_{u^{(2)}}(r)>C_{u^{(1)}}(r)$ for any $r$. At $r=0$, $C_{u^{(2)}}(0)>C_{u}(0)$ and 
$u^{(2)}(0)=u(0)$. Thus, there exists $\delta>0$ such that $u^{(2)}(r)>u(r)$ for $0<r<\delta$. We assume that 
$\delta$ is the least upper bound of such values. By contradiction, suppose that $\delta<a$. As $u^{(2)}(\delta)=u(\delta)$
and $u^{(2)'}(\delta)\leq u'(\delta)$, $\psi^{(2)}(\delta)\leq\psi(\delta)$ and 
\begin{equation}\label{alfa}
\int_0^{\delta}\frac{d}{dr}\left(\sin\psi(r)-\sin\psi^{(2)}(r)\right)dr=\sin\psi(\delta)-\sin\psi^{(2)}(\delta):=\alpha(\delta)\geq 0.
\end{equation}
Then there exists $\bar{r}\in (0,\delta)$ such that 
$$C_u(\bar{r})=(\sin\psi)'(\bar{r})>(\sin\psi^{(2)})'(\bar{r})=C_{u^{(2)}}(\bar{r}).$$
As $C_u(r)$ is increasing, $C_{u}(r)>C_{u^{(2)}}(r)$ for $r\in (\bar{r},a)$. In particular, and 
using  $u'(a)=u^{(2)'}(a)$, 
$$0<\int_{\delta}^a(C_{u}(r)-C_{u^{(2)}}(r))dr=\int_{\delta}^a\frac{d}{dr}\left(\sin\psi(r)-\sin\psi^{(2)}(r)\right)dr=-\alpha(\delta)$$
in contradiction with (\ref{alfa}).  
\end{proof}

As conclusion, the circular arcs $\Sigma^{(1)}$, $\Sigma^{(2)}$ lie respectively below and above the solution curve. From the Claim, and (\ref{u1})-(\ref{u2}),  we obtain
$$u_0+\frac{1}{\kappa u_0}-\sqrt{\frac{1}{\kappa^2 u_0^2}-a^2}<u(a)<u_0+\frac{a}{\cos\gamma}(1-\sin\gamma).$$
Together with $u_0<1/(a\kappa)$ and (\ref{q1}), we conclude

\begin{theorem} Let $S$ be a $\kappa$-cylindrical surface where $\gamma$ denotes the 
contact angle with the vertical plates $P_i$ and $0\leq\gamma<\pi/2$. Then  the difference value $q=u(a)-u(0)$  satisfies
\begin{equation}\label{q2}
\frac{1}{\kappa u_0}\left(1-\sqrt{1-a^2\kappa^2 u_0^2}\right)<q<\frac{a}{\cos\gamma}(1-\sin\gamma).
\end{equation}
\begin{equation}\label{q3}
q>\frac{2a(1-\sin\gamma)}{1+\sqrt{1+2\kappa a^2(1-\sin\gamma)}}.
\end{equation}
\end{theorem}

We can compare the upper bound for $q$ in (\ref{q2}) with
the one obtained in (\ref{q1}). So, in the case that $S$ is vertical at the walls, 
$$\frac{2a}{1+\sqrt{1+2\kappa a^2}}<q<\{a,\sqrt{\frac{2}{\kappa}}\}.$$
Another consequence of the Claim is that it allows to 
compare the volume per unit of length of $\partial\Omega$ between 
the cylindrical capillary channel and the sections of horizontal round cylinder determined
by $\Sigma^{(1)}$ and $\Sigma^{(2)}$. For this, it suffices with
\begin{equation}\label{c1}
\int_0^a  u^{(1)}(r)\ dr<\int_0^a u(r)\ dr<\int_0^a u^{(2)}(r)\ dr.
\end{equation}
The integral for $u$ can be computed by (\ref{capi3}): 
\begin{equation}\label{c2}
\int_0^a \kappa u(r)\ dr =\cos\gamma.
\end{equation}
If we denote
$$F(u_0;R)=a(R+u_0)-\frac{a}{2}\sqrt{R^2-a^2}-\frac{R^2}{2}\arcsin(\frac{a}{R}),$$
then (\ref{c1}) and (\ref{c2}) imply
$$F(u_0;R_1)<\frac{\cos\gamma}{\kappa}<F(u_0,R_2).$$
Thus, each one of the two above inequalities gives
\begin{equation}\label{A}
a\left(\frac{1}{\kappa u_0}+ u_0\right)-\frac{a}{2}\sqrt{\frac{1}{\kappa^2 u_0^2}-a^2}-\frac{\arcsin(a\kappa u_0)}{2 \kappa^2 u_0^2}<\frac{\cos\gamma}{\kappa}
\end{equation}
\begin{equation}\label{B}
\frac{\cos\gamma}{\kappa}<\frac{a^2}{\cos\gamma}+au_0-\frac{a^2\tan\gamma}{2}-\frac{a^2}{2\cos^2\gamma}\left(\frac{\pi}{2}-\gamma\right).
\end{equation}
From (\ref{B}), we obtain a lower bound for $u_0$. On the other hand, 
and since  $\partial F/\partial u_0>0$, let $u_0^+>u_0$ be the number such that $F(u_0^+;R_1)=\cos\gamma/\kappa$. As $F(x;R)-a x$ is positive, 
$$F(\frac{\cos\gamma}{a\kappa};R_1)>\frac{\cos\gamma}{\kappa}=F(u_0^+;R_1),$$
and thus 
$$u_0^+<\frac{\cos\gamma}{a\kappa}.$$

\begin{theorem}\label{tlaplace}  Let $S$ be a $\kappa$-cylindrical surface, $\kappa>0$, given 
by the profile $u=u(r;u_0)$. If $0\leq\gamma<\pi/2$ denotes the 
contact angle with the vertical plates $P_i$ at $r=a$, then 
\begin{equation}\label{laplace}
\frac{\cos\gamma}{a\kappa}-\frac{a}{\cos\gamma}+\frac{a\tan\gamma}{2}+\frac{a}{2\cos^2\gamma}\left(\frac{\pi}{2}-\gamma\right)<u_0<u_0^+<\frac{\cos\gamma}{a\kappa}.
\end{equation}
\end{theorem}

The left inequality in (\ref{laplace}) extends the one  obtained by Laplace for the 
circular capillary tube \cite{la}. The inequality $u_0<\cos\gamma/(a\kappa)$ is also a consequence 
by comparing  the slopes of $u^{(1)}$ and $u$ at the point $r=a$: $u^{(1)'}(a)<u'(a)$.
 On the other hand, the 
combination of inequalities  (\ref{A}) and (\ref{B}) gives an estimate of $u_0$ that it  
is rather cumbersome, even in the case  $\gamma=0$:
$$a\left(\frac{1}{\kappa u_0}+ u_0\right)-\frac{a}{2}\sqrt{\frac{1}{\kappa^2 u_0^2}-a^2}-\frac{\arcsin(a\kappa u_0)}{2 \kappa^2 u_0^2}<a^2+au_0-\frac{\pi a^2}{4}.$$

Now, we bound the outer height $u(a)$.
Let us move down the circular arc $\Sigma^{(2)}$ until it meets the solution curve (tangentially) at $(a,u(a))$.

{\it Claim.} At the contact point $(a,u(a))$, the arc $\Sigma^{(2)}$ lies below the solution curve $u$.

\begin{proof}[of de Claim]
The argument is similar as in the above Claim. In the new position, we compare the curvatures of 
$u$ and $\Sigma^{(2)}$: by (\ref{B2}), we have
$$C_u(a)=\kappa u(a)>\frac{\sin\psi(a)}{a}=C_{u^{(2)}}(a).$$
Thus, around the point $r=a$, $u>u^{(2)}$. By contradiction, assume that there is $\delta\in (0,a)$ such that 
$u^{(2)}(r)<u(r)$ for $r\in(\delta,a)$ and $u^{(2)}(\delta)=u(\delta)$. Since $u'(\delta)\geq u^{(2)'}(\delta)$, then
$\psi^{(2)}(\delta)\leq \psi(\delta)$. This implies 
\begin{equation}\label{C}
\int_{\delta}^a(C_{u^{(2)}}(r)-C_u(r))dr=\sin\psi(\delta)-\sin\psi^{(2)}(\delta)\geq 0.
\end{equation}
Then there would be $\bar{r}\in (\delta,a)$ such that $C_{u^{(2)}}(\bar{r})-C_u(\bar{r})>0$.
As $C_u(r)$ is increasing on $r$, $C_u(r)<C_{u^{(2)}}(r)$ on $(0,\bar{r})$ and hence also 
throughout $(0,\delta)\subset(0,\bar{r})$. Thus
$$0<\int_0 ^\delta (C_{u^{(2)}}(r)-C_u(r))dr=\sin\psi(\delta)-\sin\psi^{(2)}(\delta)\leq 0$$
by (\ref{C}). This contradiction shows the Claim. 

\end{proof}
Let $u^{(3)}$ be the displaced arc $\Sigma^{(2)}$. Then  the Claim allows to estimate the value $u(a)$ by
$$\int_0^a u^{(3)}(r) dr<\int_0^a u(r) dr.$$
Recall that the center of $u^{(3)}$ is $u_0-(u^{(2)}(a)-u(a))$. Then
$$F(u_0+u(a)-u^{(2)}(a);R_2)<\frac{\cos\gamma}{\kappa}.$$

\begin{theorem} With the same notation as in Theorem \ref{tlaplace}, for any $0\leq\gamma<\pi/2$ and $\kappa>0$, 
\begin{equation}\label{ua}
u(a)<\frac{\cos\gamma}{\kappa a}-\frac{a}{2}\tan\gamma+\frac{a}{2\cos^2\gamma}\left(\frac{\pi}{2}-\gamma\right).
\end{equation}
\end{theorem}

\begin{corollary} For any $r\in (0,a)$ and $0\leq\gamma<\pi/2$, we have
\begin{equation}\label{meniscus1}
\frac{r^2\kappa u_0}{1+\sqrt{1-r^2\kappa^2 u_0^2}}<u(r)-u_0<\frac{a}{\cos\gamma}-
\sqrt{\frac{a^2}{\cos^2\gamma}-r^2}
\end{equation}
and
\begin{equation}\label{meniscus2}
u(a)+\frac{a\sin\gamma}{\cos\gamma}-\sqrt{\frac{a^2}{\cos^2\gamma}-r^2}
<u(r)-u_0.
\end{equation}
\end{corollary}

\begin{proof} The bounds (\ref{meniscus1}) are consequence of 
$u^{(1)}(r)<u(r)<u^{(2)}(r)$. The lower bound (\ref{meniscus2}) comes from $u^{(3)}(r)<u(r)$ in $(0,a)$. 

\end{proof}

This section ends by obtaining lower estimates for the values $u_0$ and the 
difference value $q=u(a)-u(0)$. 
Since $u'>0$ in the interval $(0,a)$,  we may introduce the inclination angle $\psi=
\arctan u'(r)$ as independent variable. We have then
\begin{equation}\label{angle}
\frac{dr}{d\psi}=\frac{\cos\psi}{\kappa u}\hspace*{1cm}\frac{du}{d\psi}=\frac{\sin\psi}{\kappa u}.
\end{equation}
Simple quadratures then yield again
\begin{equation}\label{upsi}
u(\psi)=\sqrt{u_0^2+\frac{2}{\kappa}(1-\cos\psi)}
\end{equation}
obtained in (\ref{cotaz}). As a consequence, {\it the difference of squares of the maximum and minimum heights  satisfies
\begin{equation}\label{square}
u^2(\psi)-u_0^2=\frac{2}{\kappa}(1-\cos\psi).
\end{equation}
and thus, independent of the width of the strip $\Omega$}.

As $r\kappa u_0<\sin\psi$, 

\begin{corollary} In the range $0<\psi\leq\pi/2$ there holds 
$$u(\psi)<\sqrt{\left(\frac{\sin\psi}{\kappa r}\right)^2+\frac{2}{\kappa}(1-\cos\psi)}.$$
\end{corollary}
Now, the following computations are similar to the case that $\Omega$ is a 
circular disc \cite{fi0}. Let 
$$m=\cos(\psi/2),\hspace*{1cm}p=\sqrt{1+\kappa(r/m)^2}.$$
The function $r/m$ is increasing in $\psi$. As
$$u<\frac{\sin\psi}{\kappa r}p,$$
it follows from (\ref{angle}) that $p d r>r \cot\psi d\psi$, that is 
\begin{equation}\label{cot}
\frac{\sqrt{m^2+\kappa r^2}}{m r} dr >\cot\psi d\psi.
\end{equation}
From (\ref{B2}), 
$$\lim_{\psi\rightarrow 0}\frac{r(\psi)}{\sin\psi}=\frac{1}{\kappa u_0}.$$
An integration in (\ref{cot}) leads to 

\begin{theorem} \label{tt-table} In the range $0<\psi\leq\pi/2$ there holds
\begin{equation}\label{u0}
u_0>\frac{\sin\psi}{2\kappa r} \frac{\kappa}{m}(1+p)e^{1-p}.
\end{equation}
\end{theorem}

\begin{theorem}\label{t-table} There holds always for any $0\leq\gamma<\pi/2$
$$\frac{2(1-\sin\gamma)}{\kappa f(\gamma)}<u(a)-u_0<
\frac{a(1-\sin\gamma)}{\cos\gamma},$$
where
$$f(\gamma)=\frac{2\cos\gamma}{\kappa a}-
\frac{a\tan\gamma}{2}+\frac{a}{2\cos^2\gamma}\left(\frac{\pi}{2}-\gamma\right).$$
\end{theorem}

\begin{proof} The right inequality is a consequence of (\ref{q2}). For the left one, we know 
from (\ref{square}) that 
$$u(a)-u_0=\frac{2(1-\sin\gamma)}{\kappa (u(a)+u_0)}.$$
Then we bound $u(a)$ and $u_0$ by (\ref{laplace}) and (\ref{ua}). 

\end{proof}

%%%%%%%%%%%%%%%%%%%%%%%%%%%%
\section{Sessile liquid channels}\label{estimates2}
%%%%%%%%%%%%%%%%%%%%%%%%%%%%%%%

In this  section we study the setting of  a 
liquid deposited over the strip $\Omega$. As
$\kappa>0$, the vertical gravity fields points towards $\Pi$. 
We know that  $u=u(r,y)$, $(r,y)\in\Omega$, satisfies $\mbox{div } 
Tu=\kappa u$.  We assume $u(r)=u(r,y)$ and $u$ satisfies the equation 
(\ref{capi3}) with  $u(0)=u_0>0$. We write (\ref{capi3}) in terms of the 
inclination angle $\psi$ with respect to the $r$-axis:
\begin{equation}\label{con}
\frac{dr}{d\psi}=\frac{\cos\psi}{\kappa u}\hspace*{1cm}\frac{du}{d\psi}=\frac{\sin\psi}{\kappa u}.
\end{equation}
We point out that  this angle $\psi$ agrees at the 
contact between $S$ and $\Pi$ with the value $\gamma$, the angle which the 
surface meets $\Pi$ along the boundary. We know from  Section \ref{planes} that  there exists $R>0$ where $u$ is vertical, that is, 
$r(\pi/2)=R$.  Theorem \ref{t-sessile} asserts that $\psi$ takes 
any real number and the solution $u(\psi)$
 can be continued as a solution of (\ref{capi3}).

\begin{theorem} The functions $r(\psi)$ and $u(\psi)$ can be continued throughout the range $0<\psi<\pi$. Moreover, there exists  
a value $r_o=\lim_{\psi\rightarrow\pi}r(\psi)$, being $r_o>0$ and
the functions $r(\psi)$ and  $u(\psi)$ are monotonely increasing in $(0,\pi)$.
\end{theorem}

\begin{proof}
From (\ref{B2}), we know that $r<1/(\kappa u_0)$. Denote by (--) and (+) the part of the meniscus defined for 
$\psi\in (0,\pi/2)$ and $\psi\in (\pi/2,\pi)$ respectively. 
 For $r<R$ and close to $R$, we have from (\ref{B1}):
$$1-\sin\psi^-(r)=\kappa\int_r ^R u^- dt;\hspace*{1cm}1-\sin\psi^+(r)=\kappa\int_r^R u^+ dt.$$
Subtracting, we get
\begin{equation}\label{sus}
\sin\psi^-(r) -\sin\psi^+(r)=\kappa\int_r^R(u^+(t)-u^-(t))\ dt.
\end{equation}
In particular, for $r$ close to $R$, $u^+>u^-$ and hence $\sin\psi^->\sin\psi^+$. From (\ref{con}), $u^+, u^-$ are both 
increasing in $\psi$  as one can continue the solution until 
$\psi=\pi$.  However it is not possible to arrive until $r=0$, since
(\ref{sus}) would imply $0>-\sin\psi^+(0)=\kappa\int_0^R (u^+(t)-u^-(t)) dt>0$. Thus,  
$r_o=\lim_{r\rightarrow\pi}r(\psi)>0$. 

\end{proof}

The  above computation leads to

\begin{corollary} Given $\kappa, u_0>0$ and $0<\gamma\leq\pi$, there 
exists exactly one $\kappa$-cylindrical surface given by the profile $u(r;u_0)$ which makes a contact angle $\gamma$.
\end{corollary}

We now study the behavior of the sessile liquid channel with respect 
to the volume that encloses. Although our  channels  have  infinite volume, we can consider the density of fluid, that is, the volume 
per unit of length. If  $\Omega_b=(-a,a)\times (-b/2,b/2)\subset\Omega$, the volume 
of $S$ in $\Omega_b$ is 
$$2b\left(ru(a) -\int_0^a u(t)\ dt\right).$$
We call the {\it volume} of $S$ as 
$$\V=2(ru(r) -\int_0^r u(t)\ dt).$$
We write $\V(\psi)$ to denote the dependence on the angle $\psi$. Using (\ref{B1}), 
\begin{equation}\label{vo2}
\V(\psi)=2\left(ru(r)-\frac{\sin\psi}{\kappa}\right).
\end{equation}

\begin{theorem}\label{ex-vo}
 Let $V>0$ and $0<\gamma\leq\pi$. There is exactly one $\kappa$-cylindrical 
surface resting on $\Pi$, $\kappa>0$, with contact angle $\gamma$ and volume $\V=V$.
\end{theorem}

\begin{proof}
The function $\V$ is continuously differentiable on $u_0$: this follows from the
standard continuous dependence theorem of the O.D.E. theory.

We first prove existence. By (\ref{B2}), $R<1/(\kappa u_0)$. Thus $R\rightarrow 0$ as $u_0\rightarrow\infty$. We know from Section \ref{estimates1} that the function 
$u^{(3)}$ lies below $u$. Hence that a circle of radius $R$ contains 
the function $u(\psi)$, for $0<\psi<\pi/2$. Since $\sin\psi^+<\sin\psi^-$, the same circle contains also the upper part of $u$, that is, $u(\psi)$, for $\pi/2<\psi\leq\pi$.
Thus $\V<\pi R^2$, which goes to $0$ as $u_0\rightarrow\infty$.

Now, let us $u_0\rightarrow 0$. From (\ref{u0}), $r(\gamma;u_0)\rightarrow\infty$ 
for any fixed $0\leq\gamma\leq\pi/2$. By (\ref{upsi}), 
$u(\gamma;u_0)>\sqrt{2/\kappa(1-\cos\gamma)}$. Since the surface is convex, 
$\V\rightarrow\infty$ as $u_0\rightarrow 0$. In the case $\gamma>\pi/2$, 
$\V(\gamma;u_0)>\V(\pi/2;u_0)$ and the same conclusion holds.

Now, let us fix $\gamma$ and $\V$. By letting $u_0$ between $0$ and $\infty $, 
the  volume function $\V$ takes 
all values. The continuity of $\V$ with respect to $u_0$ gives the 
existence of a $\kappa$-cylindrical surface with prescribed volume $\V$.

The proof of uniqueness is obtained if we prove 
$$\stackrel{\cdot}{\V}:=\frac{\partial \V(\psi;u_0)}{\partial u_0}<0$$
for all $u_0>0$ and each fixed $\psi$ in $0<\psi\leq\pi$.
From (\ref{vo2}), 
\begin{equation}\label{volume}
\stackrel{\cdot}{\V}=2(\stackrel{\cdot}{r}u+r\stackrel{\cdot}{u}).
\end{equation}
It is known that 
$$\stackrel{\cdot}{r}(0)=0,\hspace*{1cm}\stackrel{\cdot}{u}(0)=1.$$
First, we prove the following

{\it Claim.}  
\begin{equation}\label{vpsi}
\frac{d\stackrel{\cdot}{\V}}{d\psi}<0\hspace*{.5cm}\mbox{in }(0,\pi].
\end{equation}
\begin{proof}[of the Claim]
By using (\ref{con}), 
\begin{equation}\label{vo}
\frac{d\stackrel{\cdot}{\V}}{d\psi}=\frac{2\sin\psi}{\kappa u^2}(\stackrel{\cdot}{r}u-
r\stackrel{\cdot}{u}).
\end{equation}
We shall prove that $\stackrel{\cdot}{r}<0$ and $\stackrel{\cdot}{u}>0$. 
Again,  (\ref{con}) yields
\begin{equation}\label{rupsi}
\frac{d\stackrel{\cdot}{r}}{d\psi}=-\frac{\cos\psi \stackrel{\cdot}{u}}{\kappa u^2},
\hspace*{1cm}\frac{d\stackrel{\cdot}{u}}{d\psi}=-\frac{\sin\psi \stackrel{\cdot}{u}}{\kappa u^2}.
\end{equation}
As $u_0>0$, $\frac{dr}{d\psi}(0)=-1/(\kappa u_0^2)<0$, $\stackrel{\cdot}{r}<0$ in 
an initial interval $J=(0,\delta)$, with $\delta\leq \pi$.

On the other hand, $\stackrel{\cdot}{u}>0$ for sufficiently small $\psi$ and thus, 
$$\frac{d\stackrel{\cdot}{u}}{d\psi}>-\frac{\sin\psi \stackrel{\cdot}{u}}{\kappa u_0^2}.$$
By integrating this expression, we have
\begin{equation}\label{uu}
\stackrel{\cdot}{u}>\mbox{exp }\left\{\frac{\cos\psi-1}{\kappa u_0}\right\}.
\end{equation}
We conclude $\stackrel{\cdot}{u}>0$ and (\ref{uu}) holds in $J$. From
(\ref{vo}),
\begin{equation}\label{vo1}
\frac{d \stackrel{\cdot}{\V}}{d\psi}<0\hspace*{.5cm}\mbox{in }J.
\end{equation}
By contradiction, we suppose there exists $0<\psi_0<\pi$ such that 
$\stackrel{\cdot}{r}(\psi_0)=0$. Take $\psi_0$ the first $\psi$ with
this property. As $\stackrel{\cdot}{\V}(0)=0$, (\ref{vo1}) implies
that $\stackrel{\cdot}{\V}(\psi_0)<0$. Moreover, (\ref{uu})
gives $\stackrel{\cdot}{u}(\psi_0)>0$. Now (\ref{volume}) yields
$$\stackrel{\cdot}{\V}(\psi_0)=2r(\psi_0)\stackrel{\cdot}{u}(\psi_0)>0.$$
This contradiction implies that $\stackrel{\cdot}{r}$ is negative in $(0,\pi)$and the Claim is showed. 
\end{proof}

Hence $\stackrel{\cdot}{\V}<0$ in $(0,\pi)$ for any $u_0$. 
This proves the uniqueness of Theorem. 

\end{proof}

We establish a relation between the volume $\V$ enclosed by a liquid channel with the width of the
strip $\Omega$ that defines in $\Pi$ and the contact angle $\gamma$.
 If $0<\gamma\leq\pi$ is the 
contact angle, let us denote $a=r(\gamma)$. The formulas that 
we will obtain are a consequence to compare the 
liquid channel with the halfcylinders defined by the function
$u^{(3)}$ in Section \ref{estimates1}.

We know  that the function $u^{(3)}$ is tangent 
to $u$ at the point $(a, u(a))$, $0\leq \gamma<\pi/2$. We prove that the
semicircle determined by $u^{(3)}$ in the halfplane $r>0$ contains inside
the solution curve $u$, with a single point of contact, namely, $(a,u(a))$. 
At this point, we compare the curvatures of the curves $u$ and $u^{(3)}$. By 
(\ref{B2}), $C_u(a)>C_{u^{(3)}}$. Moreover, $C_u$ is increasing on $\psi$ for 
any $0<\psi<\pi$:
$$\frac{d C_u}{d\psi}=\kappa\frac{du}{d\psi}=\frac{\sin\psi}{u}>0$$
by (\ref{con}). This proves the inclusion property.  As a consequence, we can compare the volume $\V$ of the 
liquid channel with respect to the halfcylinder determined by $u^{(3)}$. 
Denote $2R>0$ the maximal width of the liquid channel, that is, where 
the fluid is vertical at the walls.

\begin{theorem} \label{ttt-table} Let $S$ be a $\kappa$-cylindrical surface resting 
on a horizontal plane $\Pi$, $\kappa>0$, and let $\gamma$ be the angle of contact. Denote $V(\gamma)$ the enclosed volume by $S$. In the range $0<\gamma\leq\pi/2$, there holds:
\begin{equation}\label{v1}
\V(\gamma)<\frac{a^2}{\sin^2\gamma}\left(\gamma-\sin\gamma\cos\gamma\right).
\end{equation}
If $\pi/2\leq\gamma\leq\pi$, there holds
\begin{equation}\label{v2}
\V(\gamma)<R^2(\gamma-\sin\gamma\cos\gamma).
\end{equation}
\end{theorem}

\begin{proof}
For the  case $0<\gamma\leq\pi/2$, we only point that 
$$\int_0^a u^{(3)}(r)\ dr=a^2\cot\gamma+a u(a)-\frac{a^2}{\sin^2\gamma}
\left(\frac{\gamma}{2}+\frac12\sin\gamma\cos\gamma\right).$$

For (\ref{v2}), we consider a halfcircle $v$ centered at $(0,u(R))$ of radius $R$. It is known
that the lower part of this circle lies below $u$. We parametrize $v$ by 
the angle with the $r$-axis in each point. We prove that $v(\gamma)>u(\gamma)$. 
Fixed $r<R$ and $\pi/2\leq\psi\leq\pi$ with $u(\psi)=u(r)$, the function $v$ lies above 
$u$ at $r$,  $v(r)>u(r)$, and, $\sin\psi^+(r)<\sin\phi^+(r)$, where $\phi^+$ is the inclination angle of the curve $v$. As $\sin\phi^+$ decreases as $\phi^+\rightarrow
\pi$, $v(\gamma)>u(\gamma)$. Then (\ref{v2}) is a consequence of 
the computation of the volume of $v$ until $\psi=\gamma$. 

\end{proof}

We see a lower bound of the volume. 

\begin{theorem} With the same notation as in Theorem \ref{ttt-table},  there holds in the range $0<\gamma\leq\pi/2$
\begin{equation}\label{v3}
\V(\gamma)>\frac{\gamma-\sin\gamma\cos\gamma}{\kappa^2 u(\gamma)^2}.
\end{equation}
\end{theorem}

\begin{proof}
Consider the circle 
$$v(r)=u_0+R-\sqrt{R^2-r^2},\hspace*{1cm}R=\frac{1}{\kappa u(\gamma)}.$$
The curve $v$ touches tangentially $u$ at $(0,u_0)$. 
As $C_u(0)=\kappa u_0<C_v(0)=\kappa u(\gamma)$, 
$v(r)>u(r)$ for each $r$ where $v$ is defined. If we prove
that $v(\gamma)<u(\gamma)$, then the arc $v$ until $\phi=\gamma$ lies 
above $u$ and this allows to obtain a lower bound for the volume of $u$. 

At the point where $v$ attains the inclination angle $\psi=\gamma$, 
$$v(\gamma)=u_0+\frac{1-\cos\gamma}{\kappa u(\gamma)}.$$
Then $u(\gamma)>v(\gamma)$ if $u(\gamma)-u_0>v(\gamma)-u_0$. 
By using (\ref{upsi}),  we have to prove
$$\frac{2(1-\cos\gamma)/\kappa}{u_0+\sqrt{u_0^2+\frac{2}{\kappa}(1-\cos\gamma)}}>
\frac{1-\cos\gamma}{\kappa u(\gamma)},$$
or equivalently, 
$$2 u(\gamma)>u_0+\sqrt{u_0^2+\frac{2}{\kappa}(1-\cos\gamma)}.$$
But the second summand on the right side  is exactly 
$u(\gamma)$, using (\ref{upsi}) again. 
This proves that $v(\gamma)<u(\gamma)$. Then (\ref{v3}) is a consequence 
of the computation of the volume enclosed by the function $v$. 

\end{proof}

Now we prove the following inclusion result.

\begin{theorem}\label{inclusion}
Let $0<\gamma\leq\pi/2$ and let $S$ be a $\kappa$-cylindrical surface supported
on the horizontal plane $\Pi$, $\kappa>0$, and 
making contact angle $\gamma$. Let $\V$ be its volume. Then 
every $\kappa$-cylindrical surface resting on $\Pi$ and
 with smaller volume and making the same
contact angle can be translated rigidly so that it lies strictly interior 
into $S$.
\end{theorem}

\begin{proof} 
Assume that $S$ is given by the solution $u(r;u_0)$, $u_0>0$.
 Consider the solution $u^{\delta}=u(r;u_0+\delta)$, with $\delta>0$. From (\ref{B1}), 
$$\sin\psi^{\delta}-\sin\psi=\kappa\int_0^r(u^{\delta}-u)\ dt.$$
As $(u^{\delta}-u)(0)=\delta>0$, $\sin\psi^{\delta}>\sin\psi$. By (\ref{B2}) and 
(\ref{B1}), $(u^{\delta}-u)'>0$. It follows that if 
we move downward the curve $u^{\delta}$ a distance $\delta$, then it
lies above the curve $u$ except at the single point $(0,u_0)$.

The result is proved if for any $\delta>0$, $u^{\delta}<u+\delta$
at the points where the angle $\gamma$ is achieved. For given $\gamma$, 
\begin{equation}\label{upunto}
(u^{\delta}-u-\delta)(\gamma)=\int_0^{\delta}
(\stackrel{\cdot}{u}-1)\ du_0,
\end{equation}
with 
$$\stackrel{\cdot}{u}=\frac{\partial}{\partial u_0}u(\gamma;u_0).$$
Since $\stackrel{\cdot}{u}(0;u_0)=1$, 
$$\stackrel{\cdot}{u}=\int_0^{\gamma}\frac{d\stackrel{\cdot}{u}}{d\psi}\ d\psi+1.$$
As we have seen in the proof of Theorem \ref{ex-vo}), 
$\stackrel{\cdot}{u}<0$ and  
$$\frac{d\stackrel{\cdot}{u}}{d\psi}<0\hspace*{1cm}0<\psi<\gamma.$$
Thus $\stackrel{\cdot}{u}-1<0$, which implies that the integrand in (\ref{upunto})
is negative, proving the result. 

\end{proof}

We end the section obtaining  new estimates  of a
sessile  liquid
channel, with special attention if the contact angle 
lies in the range $[\pi/2,\pi]$.

\begin{theorem} \label{lipo} Let $S$ be a $\kappa$-cylindrical surface supported on
$\Pi$ and $\kappa>0$. Suppose that $u=u(\psi)$ is the profile of $S$, where 
$\psi$ denote the inclination angle with respect to the $r$-axis. 
Then in the range $0<\psi\leq\pi$ there holds
$$u(\psi)-u_0<\sqrt{\frac{2(1-\cos\psi)}{\kappa}}.$$
In the range, $\pi/2\leq\psi\leq\pi$, 
$$R-r(\psi)<\frac{1}{\sqrt{\kappa}}
\left(\sqrt{2}+\log{\left(\tan\frac{\pi}{8}\right)}\right)-2\cos{\frac{\psi}{2}}-\log{\left(\tan\frac{\psi}{4}\right)}.$$
In particular, 
$$u(\psi)-u(R)<\frac{\sqrt{2(1-\cos\psi)}-\sqrt{2}}{\sqrt{\kappa}},\hspace*{1cm}\pi/2\leq\psi\leq\pi.$$
$$R-r_o<\sqrt{\frac{2}{\kappa}},\hspace*{.5cm}u(\pi)-u(R)<\frac{2-\sqrt{2}}{\sqrt{\kappa}}.$$

\end{theorem}
\begin{proof}

By using (\ref{upsi}), we estimate  $u(\psi)$ from below as 
$$u(\psi)>\sqrt{\frac{2}{\kappa}(1-\cos\psi)}.$$
In combination with (\ref{con}), we obtain, 
$$\frac{du}{d\psi}<\frac{\sin\psi}{\sqrt{2\kappa(1-\cos\psi)}},$$
and for $\pi/2\leq\psi \leq\pi$, 
$$\frac{dr}{d\psi}>\frac{\cos\psi}{\sqrt{2\kappa(1-\cos\psi)}}.$$
The proof finishes by integrating the two above inequalities. 
\end{proof}

The bound $R-r(\psi)$ gives the minimum distance for two liquid channels in fixed parallel strips can be without  contact. One can imagine  that if  the amount of liquid is small, the shapes adopted by the 
liquid channels are graphs. If we increase the volume of 
fluid, the interfaces leave to be graphs and $\gamma>\pi/2$. Then there 
exists a critical angle where the  channels
touch their self.   Theorem \ref{lipo} gives an estimate of the distance between
each two consecutive hydrophilic strips. Other estimate is the following 

\begin{theorem} Let $S$ be a $\kappa$-cylindrical surface supported on the plane $\Pi$ and  $\kappa>0$. Assume that the contact angle $\gamma$ satisfies
$\pi/2\leq\gamma\leq\pi$. Then
$$\frac{1}{\sqrt{\kappa}}\frac{1-\sin\gamma}{\sqrt{2(1-\cos\gamma)+\kappa u_0^2}}<
R-r(\gamma)<\frac{1}{\sqrt{\kappa}}\frac{1-\sin\gamma}{\sqrt{2+\kappa u_0^2}}$$
\end{theorem}

\begin{proof} Recall that the angle parameter $\psi$ agree with the 
real contact angle $\gamma$ with $\Pi$.  As $\cos\gamma<\cos\psi<0$,
 (\ref{upsi}) gives
$$\sqrt{u_0^2+\frac{2}{\kappa}}<u(\psi)<\sqrt{u_0^2+2(1-\cos\gamma)/\kappa}.$$
Substituting into (\ref{con}), we obtain  
$$\frac{1}{\sqrt{\kappa}}\frac{\cos\psi}{\sqrt{2+\kappa u_0^2}}<
\frac{dr}{d\psi}<
\frac{1}{\sqrt{\kappa}}\frac{\cos\psi}{\sqrt{2(1-\cos\gamma)+\kappa u_0^2}},$$
and the result follows by integrating from $\psi=\pi/2$ until $\psi=\gamma$.  
\end{proof}

We can compare with the situation of absence of gravity and  pieces of  infinite cylinders,  
whose boundary is  $\partial\Omega
=L_1\cup L_2$. Then for  $\pi/2\leq\gamma\leq\pi$, 
 the amount $R-r(\gamma)$ is exactly $(1-\sin\gamma)/(2H)$, where $H$ is 
the mean curvature of the cylinder.

%%%%%%%%%%%%%%%%%%%%%%%%%%%%%%%%%%%%%%%%%%%%%%%%%%%
\section{Pendent liquid channels}\label{estimates3}
%%%%%%%%%%%%%%%%%%%%%%%%%%%%%%%%%%%%%%%%%%%%%%%%%%%

This section is devoted to the study of $\kappa$-cylindrical surfaces  when $\kappa<0$.
 In Section \ref{planes}, Theorem \ref{p-sessile},  we have studied 
its behavior.  Let 
$\alpha=(x,z)$ the directrix of the surface and without loss of generality, we assume 
$z_0=z(0)<0$. We identify $u(r(s);u_0)=z(s)$, where $u$ is a solution of 
(\ref{capi3}) with $u_0=z_0$. 
We ask when $S$ is a graph on $\Pi$, that is, if $\alpha$ is a graph on  
the $r$-axis. Theorem \ref{p-sessile} yields the necessary condition $z_0>-2/\sqrt{-\kappa}$. 
 However, in this range of values, 
it is still  possible that $S$ presents vertical points.

\begin{theorem} \label{61} Let $S$ be a $\kappa$-cylindrical surface, $\kappa<0$. Then 
$S$ is a graph on $\Pi$ if and only if 
$$-\sqrt{\frac{2}{-\kappa}}<u_0<0.$$
In such case,  there hold the following properties for the function $u=u(r;u_0)$:
\begin{enumerate}
\item The function $u$ is  periodic and it is defined on $\r$.
\item $u$ vanishes in an infinite discrete set of points.
\item The inflections of $u$ are their zeros.
\item $u_0\leq u(r)\leq-u_0$,   $u$ attains the values $\pm u_0$ and they are 
exactly the only 
critical points.
\end{enumerate}
\end{theorem}

\begin{proof}
From (\ref{cotaz}), 
$$\cos\psi=1-\frac{\kappa}{2}(z^2-z_0^2)\geq 1+\frac{\kappa}{2}z_0^2.$$
Therefore, $(x(s),z(s))$ has not vertical points if and only if  $z_0^2<-2/\kappa$. 
In such case,  $x'=\cos\psi>>0$ and  $x$ increases strictly to infinity. 
Let $s=x^{-1}(r)$.
 Using the notation of (\ref{perio}), let $r_T=x(4s_0)$. Then
$$u(r+r_T)=u(x(s)+x(4s_0))=u(x(s+4s_0))=z(s+4s_0)=z(s)=u(r).$$
This proves that $u$ is a periodic function.
Moreover, the derivative of $u$, $u'=\tan\psi$, is bounded, which implies that
 $u$ can be extended to $\r$.
From (\ref{capi3}), the inflections agree with the zeros of $u$. The rest of properties are a consequence of Theorem \ref{p-sessile}. 
\end{proof}

 We write (\ref{cotaz}) as
\begin{equation}\label{ur2}
u(\psi)^2-u_0^2=\frac{2}{\kappa}(1-\cos\psi).
\end{equation}
The case 
$$u_0=-\sqrt{\frac{2}{-\kappa}}$$
can be treated as  above, except that in a discrete set of points, 
$u$ is vertical. Moreover, (\ref{ur2}) implies that these vertical points are the
zeros of $u$. In general, we can estimate the initial interval where one can define a solution $u$ of (\ref{capi3})
without vertical points.

\begin{lemma}\label{le} Consider $u=u(r;u_0)$  the solution of (\ref{capi3})-(\ref{capi33}). Then $u$ can be
continued at least until the value $r=1/(\kappa u_0)$. Furthermore,  $\sin\psi<\kappa u_0 r$
\end{lemma}

\begin{proof}
Since $(\sin\psi)'=\kappa u$, the function $\sin\psi$   is strictly increasing on
$r$ whenever $u$ is negative. Then for $r>0$ close to $r=0$,  $\sin\psi=u'/\sqrt{1+u'^2}$
is positive. As conclusion, $u$ is increasing on $r$ near to $0$ and the expression
$$\sin\psi=\kappa\int_{0}^ru(t)\ dt$$
can be bounded in both sides by the values $u_0$ and $u(r)$. Then
\begin{equation}\label{B3}
\kappa u(r)<\frac{\sin\psi}{r}<\kappa u_0,
\end{equation}
and hence
$$\sin\psi<\kappa  u_0 r=1.$$
This means that $\psi<\pi/2$.  

\end{proof}
In general, for pendent liquid channels, we can say more about the vertical points.

\begin{theorem} Let $\alpha$ the directrix of a $\kappa$-cylindrical surface, $\kappa<0$,
such that the initial condition $z(0)=z_0$ satisfies
\begin{equation}\label{z00}
-\frac{2}{\sqrt{-\kappa}}<z_0<-\sqrt{\frac{-2}{\kappa}}.
\end{equation}
Then 
\begin{enumerate}
\item $\alpha$ presents exactly four vertical points in each arc of $\alpha$ determined
by its period. 
\item Each of these points lies 
in the segment of $\alpha$ between one extremum and one zero of $z$.
\item The height of the vertical points is $\pm\sqrt{z_0^2+2/\kappa}$.
\end{enumerate}
\end{theorem}

\begin{proof} By the symmetries of $\alpha$, it suffices to prove that between 
$s=0$ and the first time $s_0$ where $\alpha$ intersects the $r$-axis, there 
exists exactly one vertical point. Since $\theta'(s)=\kappa z(s)$, the 
function $\theta$ is increasing on $r$ in the interval $(0,s_0)$. We know
that $\theta$ attains the value $\theta=\pi/2$,
the first vertical point, at some point $s^*$, with 
$s^*<s_0$. By using again (\ref{ur2})  and (\ref{z00}), $\theta$ does not 
reach the value $\theta=\pi$.
Thus there exists a unique vertical point.
By (\ref{ur2}), the height at $s=s^*$ is 
$-\sqrt{z_0^2+\frac{2}{\kappa}}$.  
\end{proof}

When $\alpha$ begins from $s=0$, $\alpha$ is a graph on the $r$-axis until that 
$\alpha$ is vertical.
It is possible to determine  the region where 
this first vertical point occurs. To this end, one 
 can carry as in the case of  pendent liquid drops. We refer \cite{cf1}
and \cite[Ch. 4.6]{fi1}. The main argument is a "Comparison Lemma" that compares 
 $u$ with circular arcs and the hyperbolas
$ru<1/(2\kappa)$ and $ru<1/\kappa$. For example, one can show that
{\it the directrix $\alpha$, in the initial region $z<0$, does not enter the region 
$ru\leq 1/\kappa$}. We omit the details.

We summarize then the behavior of $u$, a solution of (\ref{capi3}), with  
$u_0>-\sqrt{-2/\kappa}$. After $r=0$, $u$ increases on $r$ until that 
it touches the $r$-axis at some point $R$. Theorem \ref{t-ii} says 
that $u$ is symmetric with respect to the point $(R,0)$. Thus, 
$u$ increases until the value $r=2R$, where $u$ takes the value $-u_0$. 
Again, the symmetry of $u$ with respect to the line $r=2R$ implies
that $u$ decreases until to arrive at $r=4R$ to reach the value 
$u_0$. From this position, the curve $u$ repeats the same behavior by the 
periodicity of $u$ (recall that the period is $R_T=4R$,  see Theorem \ref{61}).
The next result gives an estimate  of the value $r=R$ in the 
sense that, fixed the constant $\kappa$,  the first zero of $u$ remains bounded in 
some interval, independent on the initial value $u_0$.

\begin{theorem} Let $\kappa<0$ and $u(r;u_0)$ a solution of 
(\ref{capi3})-(\ref{capi33}) with
\begin{equation}\label{hypo}
-\sqrt{\frac{-2}{\kappa}}<u_0<0.
\end{equation}
Then
\begin{equation}\label{RR}
\frac{1}{\sqrt{-2\kappa}}<R<\sqrt{\frac{-2e}{\kappa}}.
\end{equation}
\end{theorem}

\begin{proof} 
Since $u(R)=0$, (\ref{ur2}) implies that $\cos\psi_R=1+\kappa u_0^2/2$.
From (\ref{B3}), 
$$R>\frac{\sin\psi(R)}{\kappa u_0}=
\frac{1}{\kappa u_0}\sqrt{1-\cos^2\psi_R}=\frac{-1}{2\kappa}\sqrt{-4\kappa-\kappa^2 u_0^2}.$$
The left side in (\ref{RR}) is then a consequence of this inequality  and (\ref{hypo}).
Now, we show the right inequality in (\ref{RR}). In the region where $u<0$, 
$\sin\psi$ is increasing on $r$. Let us fix $a$ such that 
$0<a<R$. Then $\sin\psi(r)>\sin\psi(a)$. As $u'(r)=\tan\psi$ 
and $\sin\psi<\tan\psi$, we have
$$u'=\tan\psi>\sin\psi\geq \frac{a}{r}\sin\psi(a).$$ 
A simple integration between $a$ and $R$ gives
$$R<a\exp{\left\{\frac{-u(a)}{\sin\psi(a)}\right\}}.$$
Again (\ref{B3}) leads to
\begin{equation}\label{exp}
R<a\exp{\left\{\frac{-1}{\kappa a^2}\right\}}.
\end{equation}
Since this holds for every $a<R$ and the function on $a$ on the right side of (\ref{exp})
 attains a minimum at $a=\sqrt{-2/\kappa}$, we 
obtain the desired estimate. 
\end{proof}

Following the same steps as in \cite{cf1,fi1}, one could improve the upper bound for $R$.

We analyze the case 
\begin{equation}\label{casou0}
-\frac{2}{\sqrt{-\kappa}}<z_0<-\sqrt{\frac{-2}{\kappa}}.
\end{equation}

\begin{theorem}\label{cinco} Let $\alpha=(x(s),z(s))$ be the directrix  of a
$\kappa$-cylindrical surface $S$. Assume that $z_0$ satisfies (\ref{casou0}).
Then there exist numbers $z_1, z_2$, with 
$$-\frac{2}{\sqrt{-\kappa}}<z_2<z_1<-\frac{\sqrt{2}}{\sqrt{-\kappa}}$$
 and the following properties hold:
\begin{enumerate}
\item If $z_1<z_0$, 
then $\alpha$ has not double points and $x$ goes to $\infty$. 
\item If $z_0=z_1$, then $\alpha$ has double points, where $\alpha$ tangentially 
meets itself at these points, $\alpha$ lies in
$\{x\geq 0\}$, and $x$ goes to $\infty$. 
\item If $z_2<z_0<z_1$, $\alpha$ has double points, meeting at these points transversally and 
$x$ goes to $\infty$. 
\item If $z_0=z_2$, then $\alpha$ is a closed curve with self intersection at the
origin. 
\item If $z_0<z_2$, $\alpha$ has double points, where $\alpha$ meets itself 
transversally and  $x$ goes to $-\infty$. 
\end{enumerate}
\end{theorem}
\begin{proof} Denote $r(\pi/2)$ and $r(\psi_0)$ the abcisas of the first 
vertical point and the first point which $\alpha$ meets the $r$-axis. We know that $r(\pi/2)>0$.

{\it Claim 1.} $r(\pi/2)\leq \frac{\pi}{2 \sqrt{-2\kappa}}$, independent of the value $z_0$.

\begin{proof}[of the Claim 1] By (\ref{casou0}), $z_0^2>-2/\kappa$. 
For each $0\leq\psi\leq\psi/2$ and by using  (\ref{ur2}), we obtain
$$\kappa u(\psi)>-\kappa\sqrt{\frac{-2\cos\psi}{\kappa}}.$$
Since $\cos\psi>0$, we have from (\ref{con}) that
$$\frac{dr}{d\psi}<\frac{\sqrt{\cos\psi}}{\sqrt{-2\kappa}}<\frac{1}{\sqrt{-2\kappa}}.$$
Integrating from $\psi=0$ to $\psi=\pi/2$, we show the Claim 1.  

\end{proof}

{\it Claim 2.} There exists a continuous function $\varphi=\varphi(z_0)$
strictly decreasing on $z_0$
such that $r(\psi_0(z_0))<\varphi(z_0)$, and 
$$\lim_{z_0\rightarrow -2/\sqrt{-\kappa}}\varphi(z_0)=-\infty.$$

\begin{proof}[of the Claim 2.] Consider $\psi\in[\pi/2,\psi_0]$.
Since $u(\psi_0)=0$, by (\ref{cotaz}), we have
$$1-\cos\psi_0=-\frac{\kappa}{2}z_0^2.$$
This proves that as $z_0\rightarrow-\frac{2}{\sqrt{-\kappa}}$, the angle $\psi_0$ which
the directrix $\alpha$ meets the $r$-axis goes to $\psi=\pi$. As $z_0^2<-4/\kappa$, 
again (\ref{ur2}) leads to
$$u(\psi)\geq -\sqrt{\frac{-2(1+\cos\psi)}{\kappa}}.$$
Then (\ref{con}) implies
$$\frac{dr}{d\psi}<\frac{1}{\sqrt{-2\kappa}}\frac{\cos\psi}{\sqrt{1+\cos\psi}}.$$
Integrating from $\psi=\pi/2$ until $\psi= \psi_0$, we obtain
$$r(\psi_0)-r(\pi/2)<\frac{1}{\sqrt{-2\kappa}}\int_{\pi/2}^{\psi_0}
\frac{\cos\psi}{\sqrt{1+\cos\psi}}d\psi.$$
An integration gives
$$r(\psi_0)-r(\pi/2)<\frac{2}{\sqrt{-\kappa}}\left(\sin(\psi_0/2)-
\mbox{arctanh}(\tan(\psi_0/4))-\frac{\sqrt{2}}{2}+\mbox{arctanh}(\tan(\pi/8))\right).$$
From the Claim 1, $r(\pi/2)$ is bounded. Then, up a constant $C$, 
$$r(\psi_0)<\frac{2}{\sqrt{-\kappa}}\left(\sin(\psi_0/2)-
\mbox{arctanh}(\tan(\psi_0/4))\right)+C:=\varphi(z_0).$$
Finally, because $\psi_0\rightarrow\pi$, we have
$\varphi\rightarrow-\infty$ as $z_0\rightarrow -2/\sqrt{\kappa}$. 

\end{proof}

 We know from (\ref{ur2}) that
if $u_a<u_b$, then $\cos\psi_0(u_a)<\cos\psi_0(u_b)$. Moreover, 

{\it Claim 3 .} The function $r(\psi_0(z_0))$ is strictly decreasing on $z_0$. 

\begin{proof} Consider $u_a<u_b$. A reasoning similar as in Theorem \ref{inclusion}
proves that if $\delta<0$, $u(r;z_0+\delta)+\delta>u(r;z_0)$ (by the Lemma 
\ref{le}, $r(\pi/2;z_0+\delta)<r(\pi/2;z_0)$). This shows that if we move upwards
$\alpha_a$ until to arrive the point $(0,u_b)$, $\alpha_a$ lies over $\alpha_b$ 
at least until the first vertical point of $\alpha_a$. Then $\alpha_a$, in the new position, lies over $\alpha_b$ at least until that both $\alpha_a$ and $\alpha_b$ meet
the $r$-axis. If $\bar{r}$ is the $x$-coordinate  of the point where (the displaced) $\alpha_a$
intersects the $r$-axis in the first time, we have, 
$$r(\psi_0(u_a))<\bar{r}<r(\psi_0(u_b)).$$ 
\end{proof}

Now, we sketch  the proof of the Theorem and we omit the details.
Take $z_0$ varying from $z_0=-\sqrt{-2/\kappa}$ until $z_0=-2/\sqrt{-\kappa}$.
Let $z_1$ and $z_2$ be the unique numbers, $z_1<z_2$ such that, in the notation
of (\ref{perio}), there hold
$$x(2s_0(z_1))=0\hspace*{.5cm}\mbox{and}\hspace*{.5cm} r(\psi_0(z_2))=0.$$
 The existence is given by the Claim 2 and the uniqueness by the Claim 3. By (\ref{perio}), the  direction,
 left or right, that takes 
$\alpha$ depends on the sign of  $x(4s_0)$. The critical time occurs
when $x(4s_0)=x(s_0)=0=r(\psi_0(z_2))$, where $\alpha$ is a closed curve. 
Moreover, since $z$ is increasing in $(0,2s_0)$
(and $z$ goes from $z_0$ to $-z_0$), $\alpha$ does not intersect itself.
Then the results follow using
the symmetries of the directrix according the Theorems \ref{t-i} and \ref{t-ii}.

\end{proof}

\begin{remark} The results obtained here show the contrast of
behavior between pendent liquid channels and pendent liquid rotational drops. In the latter 
setting, for the values of $z_0<0$ where $\alpha$ presents vertical points,  the number of vertical points goes increasing as $z_0\rightarrow-\infty$ \cite{cf1}.
However, in our case, the periodicity of the curve $\alpha$ simplifies the scene.
\end{remark}

%%%%%%%%%%%%%%%%%%%%

\end{document}